\documentclass[final]{siamltex}
\usepackage{texdraw,amsmath,amssymb,bm,cite,graphicx,color}
\newcommand{\eq}{\begin{equation}}
\newcommand{\eeq}{\end{equation}}
\newcommand{\eqn}{\begin{eqnarray}}
\newcommand{\eeqn}{\end{eqnarray}}

\newcommand{\bsea}{\begin{subeqnarray}}
\newcommand{\esea}{\end{subeqnarray}}
\newcommand{\nn}{\nonumber}
\newcommand{\nnl}{\nonumber \\}
\newcommand{\alg}[1]{\begin{align}  #1 \end{align}}

\def\bmat{\left[ \begin{array}}
\def\emat{\end{array} \right]}

\newcommand{\tr}{\mathop{\rm tr}}  

\newcommand{\Bc}{ \mathcal{B}}

\newcommand{\Dc}{ \mathcal{D}}

\newcommand{\Gc}{ \mathcal{G}}
\newcommand{\Hc}{ \mathcal{H}}

\newcommand{\Jc}{ \mathcal{J}}
\newcommand{\Kc}{ \mathcal{K}}
\newcommand{\Lc}{ \mathcal{L}}

\newcommand{\Nc}{ \mathcal{N}}
\newcommand{\Oc}{ \mathcal{O}}

\newcommand{\Qc}{ \mathcal{Q}}
\newcommand{\Rc}{ \mathcal{R}}
\newcommand{\Sc}{ \mathcal{S}}

\newcommand{\Ds}{ \mathbb{D}}
\newcommand{\Es}{ \mathbb{E}}

\newcommand{\Rs}{ \mathbb{R}}

\def\qed{\hfill \vrule height 7pt width 7pt depth 0pt \smallskip}

\newcounter{pippo}

\newtheorem{remark}{Remark}[section]
\newtheorem{teor}{Theorem}[section]
\newtheorem{corr}{Corollary}[section]
\newtheorem{propo}{Proposition}[section]
\newtheorem{lemm}{Lemma}[section]
\newtheorem{exam}{Example}
\newtheorem{probl}[pippo]{Problem}
\newtheorem{defn}{Definition}[section]
\newcommand{\teo}{\begin{teor}}
\newcommand{\eteo}{\end{teor}}
\newcommand{\cor}{\begin{corr}}
\newcommand{\ecor}{\end{corr}}
\newcommand{\prop}{\begin{propo}}
\newcommand{\eprop}{\end{propo}}
\newcommand{\lem}{\begin{lemm}}
\newcommand{\elem}{\end{lemm}}
\newcommand{\ex}{\begin{exam}}
\newcommand{\eex}{\end{exam}}
\newcommand{\pb}{\begin{probl}}
\newcommand{\epb}{\end{probl}}
\newcommand{\df}{\begin{defn}}
\newcommand{\edf}{\end{defn}}
\newcommand{\aprop}{\begin{apropo}}
\newcommand{\eaprop}{\end{apropo}}
\newcommand{\alem}{\begin{alemm}}
\newcommand{\ealem}{\end{alemm}}
\newcommand{\rem}{\begin{remark}}
\newcommand{\erem}{\end{remark}}

\title{Convergence analysis of a family of robust Kalman filters based on the contraction principle\thanks{This work has been partially supported by the FIRB project ``Learning
meets time'' (RBFR12M3AC) funded by MIUR.}}

\author{Mattia Zorzi
 \thanks{M. Zorzi is with the
Dipartimento di Ingegneria dell'Informazione, Universit\`a di
Padova, via Gradenigo 6/B, 35131 Padova, Italy,
({email: \tt\small zorzimat@dei.unipd.it}).} }

\begin{document}

\pagestyle{myheadings}
\thispagestyle{plain}
\markboth{M. ZORZI}{CONVERGENCE OF ROBUST KALMAN FILTERS}

\maketitle

\begin{abstract} 
In this paper we analyze  the convergence  of a family of robust Kalman filters. For each filter of this family the model uncertainty is tuned according to the so called tolerance parameter. Assuming that the corresponding state-space model is  reachable and observable, we show that the corresponding Riccati-like mapping is strictly contractive provided that the tolerance is sufficiently small, accordingly the filter converges.  

\end{abstract}



 \section{Introduction}
Robust Kalman filtering is a computational tool with widespread applications in many fields, e.g. \cite{ZORZI_ZENERE_IFAC}.  In this paper we consider the parametric family of robust Kalman filters introduced in \cite{STATETAU_2017}, see also the former works \cite{ROBUST_STATE_SPACE_LEVY_NIKOUKHAH_2013},\cite{LEVY_NIKOUKHAH_2004},\cite{ROBUSTNESS_HANSENSARGENT_2008}. The parameter describing this family is denoted by $\tau$. Once $\tau$ is fixed, the model uncertainty is represented by a ball  which is about the nominal model and  formed by placing a bound on the $\tau$-divergence, \cite{OPTIMALITY_ZORZI},\cite{DUAL},\cite{OPTIMAL_PREDICTION_ZORZI_2014}, between the actual and the nominal model. The  bound is fixed by the user and represents the tolerance of the mismatch between the actual and the nominal model. Then, the robust filter is obtained by minimizing the mean square error according to the least favorable model in this ball. Interestingly, relaxing the assumption that the actual model belongs to the ball, we obtain a family of risk sensitive filters parametrized by $\tau$ wherein the tolerance parameter is replaced by the so called risk sensitivity parameter. In particular, for $\tau=0$ we obtain the usual risk sensitive filter, see \cite{boel2002robustness,OPTIMAL_SPEYER_FAN_BANAVAR_1992,RISK_WHITTLE_1980,
H_INF_HASSIBI_SAYED_KAILATH_1999}.

In this paper we analyze  the convergence of this family of discrete-time robust Kalman filters.  More precisely, we prove that 
the error covariance, obeying to a Riccati-like iteration, converges to a unique positive definite solution.

The convergence of Riccati-like iterations can be performed using classical argumentations, \cite{ferrante2013generalised,Ferrante20141176,FERRANTE1996359}.  
Alternatively, the convergence analysis can be performed using the contraction principle as in the former paper by Bougerol  \cite{BOUGEROL_1993}. More precisely, under reachability and observability assumptions, he proved that the discrete-time Riccati iteration is a strict contraction for the Riemann metric associated to the cone of positive definite matrices. Interestingly, the same result holds using the Thompson's part metric  \cite{LEE_LIM_2008,GAUBERT_2012}.
The latter metric is more effective than the former in the sense that it gives a tighter bound on the convergence rate of the iteration. It is also worth noting that the contraction principle has been used also to prove the convergence of different kinds of nonlinear iterations \cite{LAWSON_LIM_2006,LAWSON_LIM_2007,LEE_LIM_2008}.

The convergence analysis that we present here is based on the contraction principle. This analysis takes the root  from the paper \cite{LEVY_ZORZI_RISK_CONTRACTION}. The latter studies the convergence of the risk sensitive Riccati iteration corresponding to the usual risk sensitive filter. In particular, placing an upper bound on the risk sensitivity parameter it is possible to prove that the $N$-fold composition of the risk sensitive Riccati mapping is strictly contractive for the Thompson's part metric. Since the robust Kalman filter with $\tau=0$ can be understood as the usual risk sensitive filter where the risk sensitivity parameter is now time-varying, it is possible to characterize an upper bound on the tolerance of this robust filter in such a way that the time-varying risk sensitivity parameter is sufficiently small. In this way, the $N$-fold composition of the mapping is strictly contractive and thus the robust filter converges, \cite{ZORZI_CONTRACTION_CDC}. In this paper we extend these results for the entire family of robust Kalman filters.

The outline of the paper is as follows. In  Section \ref{sec_thompson_contraction} we recall the Thompson's part metric for positive definite matrices and the properties of contraction mappings. In Section \ref{sec_rob_kalman} we review the robust Kalman filter, we derive the downsampled version and the corresponding $N$-fold Riccati iteration. In this way we are able to derive a condition for which the iteration is strictly contractive. In Section
\ref{sec_toler} we translate this condition in terms of upper bound on the tolerance of the robust filter.
In Section  \ref{sec_example} an illustrative example is provided. In Section \ref{sec_RS} deals with the convergence analysis of the family of $\tau$-risk sensitive filters.  Finally, we draw the conclusions in Section \ref{sec_conclusions}.

{\em Notation}. Given $x\in\Rs^n$, $\|x\|$ denotes the Euclidean norm of $x$, and $\|x\|_K$ 
denotes the weighted Euclidean norm with weight matrix $K$ positive definite. The $i$-th singular value of $P\in\Rs^{n\times n}$
is denoted by $\sigma_i(P)$ and $\sigma_1(P)\geq \sigma_2(P)\geq \ldots \geq \sigma_n(P)$. $\|P\|$ denotes the spectral norm of $P$, i.e. $\|P \|=\sigma_1(P)$. $\Qc^n$ denotes the vector space of symmetric matrices of dimension $n \times n$. The cone of positive definite matrices in $\Qc^n$ is denoted by $\Qc_+^n$, and its closure by $\bar \Qc_+^n$. $\diag(d_1\ldots d_n)$ denotes the diagonal matrix with elements in the main diagonal $d_1,\ldots,d_n$; similarly $\mathrm{blkdiag}(P_1\ldots P_n)$ denotes the block-diagonal matrix with matrices in the main block-diagonal $P_1,\ldots,P_n$. Given $P\in\Qc_+^n$ with eigendecomposition $P=UDU^T$ such that $U$ is an orthogonal matrix and $D=\diag(\sigma_1(P)\ldots \sigma_{n}(P))$, the exponentiation of $P$ to a real number $\tau$ is defined as $P^{\tau}=UD^\tau U^T$ with $D^\tau=\diag(\sigma_1(P)^\tau\ldots \sigma_{n}(P)^\tau)$. Similarly, we define $\exp(P)=U\exp(D)U^T$ with $\exp(D)=\diag(e^{\sigma_1(P)}\ldots e^{\sigma_{n}(P)})$ and $\log (P)=U\log(D)U^T$ with $\log(D)=\diag(\log(\sigma_1(P))\ldots \log(\sigma_{n}(P)))$.

\section{Thompson's part metric and contraction mappings}\label{sec_thompson_contraction}
Let $P$ and $Q$ belong to $\Qc_+^n$. The Thompson's part metric \cite{BATHIA_2003} between $P$ and $Q$ is defined as  
 \alg{
d_T(P,Q) &= \| \log (P^{-1/2} QP^{-1/2})\| \nn\\ &=\max\{\log(\sigma_1(P^{-1}Q)),\log(\sigma_1(Q^{-1}P)) \}
.\nn }
Beside all the traditional properties of a distance, $d_T$ has the feature that it is invariant under matrix inversion and congruence transformations.

Let $h(\cdot)$ be an arbitrary mapping in $\Qc_+^n$. We say that $h$ is strictly contractive if its
contraction coefficient (or {\em Lipschitz} constant)  \alg{
\xi(h) = \sup_{P,Q \in \Qc_n^+, P \neq Q} \frac{d_T(h(P),h(Q))}{d_T(P,Q)}
\nn} is less than one. Since the metric space $(\Qc_+^n,d_T)$ is complete \cite{thompson1963certain}, if $h$ is a strict contraction of $\Qc_+^n$ for the distance $d_T$, by the {\em Banach} fixed point theorem, \cite[p. 244]{AE}, 
there exists a unique fixed point $P$ of $h$ in $\Qc_+^n$ satisfying $P=h(P)$. Moreover, this fixed point is given by performing the iteration $P_{k+1}=h(P_k)$ starting with any $P_0\in\Qc_+^n$. Consider the downsampled  iteration $P_{k+1}^d=h_k^N(P_k^d)$ where $P_k^d=P_{kN}$ and $N$ is an integer. Here, $h_k^N$ is the $N$-fold composition of $h$ at step $kN$. If $h_k^N$ is strictly contractive for $k\geq \tilde q$ with $\tilde q$ fixed, then $h$ has a unique fixed point given as before. In this paper we will need the next Lemma \cite[Th. 5.3]{LEE_LIM_2008}.
\lem \label{lemma_lee_lim}Let $W_1,W_2\in\Qc_+^n$. Then, the mapping
\alg{ h(P)=M(P^{-1}+W_1)^{-1}M^T+W_2 \nn}
is strictly contractive with
\alg{\xi(h)\leq  \left(\frac{\sqrt{\sigma_1(W_1^{-1}M^T W_2^{-1} M)}}{1+\sqrt{1+\sigma_1(W_1^{-1}M^T W_2^{-1}M)}}\right)^2.\nn }
\elem It is worth noting that the results outlined in this Section also hold using the Riemann metric \cite{BOUGEROL_1993}. On the other hand, the Thompson's part metric is more effective than the Riemann one because it provides a tighter bound on the convergence rate of the previous iteration.

\section{Contraction property of the robust Kalman filters} \label{sec_rob_kalman}
Consider the state-space model 
\alg{\label{state_space_model}x_{k+1}&= A x_k+B v_k\nn\\
y_k&= Cx_k+Dv_k}
where $x_k\in\Rs^n$ is the state process, $y_k\in\Rs^p$ is the observation process and $v_k\in\Rs^m$ is white Gaussian noise with unit variance, i.e. $\Es[v_kv_k^T]=I_m$. The initial state $x_0$ is assumed to be independent of $v_k$. Moreover, its nominal probability density is $f_0(x_0)\sim \Nc(\hat x_0,V_0)$. 
Model (\ref{state_space_model}) is completely described by the nominal 
joint Gaussian probability density $ f_k(x_{k+1},y_k| Y_{k-1})$ of $x_{k+1}$ and $y_k$ conditioned on $Y_{k-1}:=[\, y_0\, \ldots y_{k-1}\,]^T$.
We consider the family of robust Kalman filters \cite{STATETAU_2017},\cite{ROBUST_STATE_SPACE_LEVY_NIKOUKHAH_2013} parametrized by $\tau\in [0,1]$: 
\alg{\hat x_{k+1}=  \underset{g_k\in\Gc_k}{\mathrm{argmin}}\underset{\tilde f_k\in   \Bc^c_{k,\tau}}{\max} \Es_{\tilde f_k}[ \|x_{k+1}-g_k(y_k)\|^2|Y_{k-1}]\nn}
where $ \Es_{\tilde f_k}[\cdot |Y_{k-1}]$ is the conditional expectation taken with respect to $\tilde f_k(x_{k+1},y_k|Y_{k-1})$ which is the least-favorable joint Gaussian probability density of $x_{k+1}$ and $y_k$ conditioned on $Y_{k-1}$.
$  \Bc^c_{k,\tau}$ is a ball about the nominal density $  f_k(x_{k+1},y_k| Y_{k-1})$ with radius $c$:
\alg{ \Bc^c_{k,\tau}=\{\tilde f_k(z_k|Y_{k-1}) \hbox{ s.t. } \Dc_\tau(\tilde f_k\|  f_k )\leq c\} \nn }
where $\Dc_\tau$ is the $\tau$-divergence family with parameter $\tau\in[0,1]$
and defined as follows. Let $\tilde f$ and $f$ be two $q$-dimensional Gaussian probability densities with mean vector $\tilde m_z,m_z$ and covariance matrix $\tilde K_z,K_z$, respectively.  Then, the $\tau$-divergence family is defined as
{\small \alg{ \mathcal{D}_\tau &(\tilde f\| f)=\nn\\
& \left\{
                               \begin{array}{ll}
                                 \|\Delta m_z\|^2_{K_Z^{-1}}+\tr\left(-\log(\tilde K_z K_z^{-1})+\tilde K_z K_z^{-1}-I_{q} \right), & \hbox{$\tau=0$} \\
                                	\frac{1}{1-\tau} \|\Delta m_z\|^2_{K_Z^{-1}}+\tr\left(\frac{1}{\tau(\tau-1)}(L_z^{-1}\tilde K_z L_z^{-T})^{\tau}  +\frac{1}{1-\tau} \tilde K_z K_z^{-1}+\frac{1}{\tau}I_{q}\right), &\hbox{$0<\tau<1$} \\
                               \delta(\Delta m_z)+  \tr\left( L_z^{-1}\tilde K_z L_z^{-T}\log(L_z^{-1}\tilde K_z                                    L_z^{-T})-\tilde K_z K_z^{-1}+I_{q}\right), & \hbox{$\tau=1$ }
                               \end{array}
                             \right.\nn }}
                             where $\Delta m_z=m_z-\tilde m_z$ and $L_z$ is such that $K_z=L_zL^T_z$.
 Note that, $\mathcal{D}_\tau (\tilde f\| f)$ coincides with the {\em Kullback-Leibler} divergence for $\tau=0$, \cite{ALPHA}.  To understand the role of parameter $\tau$ in $\mathcal{B}_{k,\tau}^c$ consider the ball $\mathcal{B}_\tau^c:=\{ \tilde f\,: \, \mathcal{D}_\tau (\tilde f\| f)\leq c	\}$. In \cite{OPTIMALITY_ZORZI}, it has been shown that, increasing $\tau$ and choosing $c$ in such a way that the measure of $\mathcal{B}_\tau^c$ remains constant, then the uncertainty described by $\mathcal{B}_\tau^c$ increases for the covariance matrix while it decreases for the mean vector. Accordingly, $\tau$ tunes how to allocate the mismodeling budget between the mean vector and the covariance matrix. $c$ is referred to as tolerance and measures the model uncertainty.
$\Gc_k$ is the class of estimators with finite second-order moments with respect to all densities $\tilde f_k(x_{k+1},y_k|Y_{k-1})\in \Bc^c_{k,\tau}$. The resulting estimator obeys the recursion: 
  \alg{\label{rob_est_formula1} \hat x_{k+1}=A \hat x_k+G_k(y_k-C \hat x_k) }
where $G_k$ is the gain matrix   
\alg{\label{def_Gk}G_k&=(A V_k C^T+B D^T)(C V_k C^T+D D^T)^{-1}.}
If $x_k-\hat x_k$ denotes the state prediction error at time $k$, its pseudo-nominal and least-favorable covariance matrix is denoted by $P_k$ and $V_k$, respectively. Then, the latter obey to the Riccati-like iteration: 
\alg{\label{def_P_t+1} P_{k+1}&= AV_k A^T- G_k(C V_k C^T+D D^T) G_k^T+B B^T \\
\label{def_V_t+1} V_{k+1}&=\left\{
      \begin{array}{ll}
        L_{P_{k+1}}\left(I_n-\theta_k(1-\tau) L_{P_{k+1}}^TL_{P_{k+1}}\right)^{\frac{1}{\tau-1}} L_{P_{k+1}}^T , & \hspace{-0.1cm}  0\leq \tau<1\\
L_{P_{k+1}}\exp\left(\theta_k L_{P_{k+1}}^T L_{P_{k+1}}\right) L_{P_{k+1}}^T, & \tau=1
      \end{array}
    \right.
}
 where $L_{P_{k+1}}$ is such that $P_{k+1}=L_{P_{k+1}}L_{P_{k+1}}^T$ and $\theta_k^{-1}>(1-\tau)\|P_{k+1}\|$
is the unique solution to \alg{ \label{cond_on_theta_t} c=\gamma_\tau(P_{k+1},\theta_k)} where
$\gamma_\tau$ is defined as: 
\alg{ \label{def_gamma} & \gamma_\tau(P,\theta)=  \left\{
      \begin{array}{ll}
       -\log\det(I_n-\theta P)^{-1}+\tr((I_n-\theta P)^{-1}-I_n) , &  \tau=0\\
\tr (-\frac{1}{\tau(1-\tau)}(I_n-\theta(1-\tau) L_{P}^TL_{P})^{\frac{\tau}{\tau-1}}  \\
 \hspace{0.3cm}+\frac{1}{1-\tau} (I_n-\theta(1-\tau) L_{P}^TL_{P})^{\frac{1}{\tau-1}}+\frac{1}{\tau}I_n)        , & \hspace{-0.6cm}0 <\tau<1\\
      \tr(\exp(\theta L_{P}^T L_{P})(\theta L_{P}^T L_{P}-I_n)+I_n)  , & \tau=1.
      \end{array}
    \right.} $\theta_k$ is called risk sensitivity parameter and it is time-varying.
    In the case that $c=0$, i.e. no uncertainty in the nominal model, we obtain the usual Kalman filter. Regarding the performance analysis of this family of robust Kalman filters with respect to parameter $\tau$ we refer to \cite{STATETAU_2017}. It is worth noting, in view of (\ref{def_V_t+1}), we have that $P_{k+1}< V_{k+1}$.  To study the asymptotic behavior of this robust Kalman filter, the matrices $A$, $B$, $C$, $D$ and the tolerance $c$ 
    are assumed to be constant. Without loss of generality we assume that $BD^T=0$. Otherwise, we can rewrite the filter (\ref{rob_est_formula1})-(\ref{def_gamma}) with $\tilde A= A-BD^T(DD^T)^{-1}C$, $\tilde B$ such that  $\tilde B\tilde B^T=B(I-D^T(DD^T)^{-1}D)B^T$, $\tilde C=C$ and $\tilde D=D$. In this way $\tilde B \tilde D^T=0$. Substituting 
 (\ref{def_Gk}) in (\ref{def_P_t+1}) and using the Woodbury formula, we obtain the Riccati-like iteration
 \alg{  \label{iteration}P_{k+1}=r_{\tau,c}(P_k):=A(V_k^{-1}+ C ^T(DD^T)^{-1} C)^{-1}A^T+B B^T.}
Defining the positive definite matrix 
\alg{  \Phi_{k}=P_{k+1}^{-1}-V_{k+1}^{-1}\nn}
we have 
\alg{ \label{map_in_Qt} r_{\tau,c}(P_k)&=A ( P_k^{-1}-\Phi_{k-1}+ C ^T(DD^T)^{-1} C)^{-1}A^T+BB^T.}

The mapping in (\ref{map_in_Qt}) has the same structure of the risk sensitive Riccati mapping, \cite{RISK_WHITTLE_1980}. Accordingly, the robust filter (\ref{rob_est_formula1})-(\ref{def_gamma}) can be interpreted as solving a standard
least-square filtering problem with time-varying parameters in {\em Krein} space, \cite{RISK_KREIN_SPACE1_HASSIBI_SAYED_KAILATH_1996,RISK_KREIN_SPACE2_HASSIBI_SAYED_KAILATH_1996}.
The {\em Krein} state-space model consists of dynamics
and observations in (\ref{state_space_model}), to which we must adjoin the
new observations
$ 0 = x_k + u_k$.
The components of noise vectors $v_k$ and $u_k$
now belong to a {\em Krein} space and have the inner product
\alg{
\left \langle \bmat {c}
v_k \\
u_k
\emat
\, , \,
\bmat {c}
v_s \\
u_s
\emat \right \rangle
= \bmat {cc}
I_m  & 0  \\
  0 & -\Phi_{k-1}^{-1}  \\
\emat \delta_{k-s} 
\nn} where $\delta_k$ denotes the Kronecker delta function. 
Since $x_k$ is
{\em Gauss-Markov}, the downsampled process $x_k^d := x_{kN}$, with
$N$ integer, is also {\em Gauss-Markov} with state-space model
\alg{
x_{k+1}^d &= A^N x_k^d + {\cal R}_N {\bf v}_k^N \nn \\
{\bf y}_k^N &= {\cal O}_N x_k^d +{\cal D}_N {\bf v}_k^N+ {\cal H}_N {\bf
v}_k^N \nn \\
{\bf 0} &= {\cal O}_N^R x_k^d + {\bf u}_k^{N} + {\cal L}_N {\bf v}_k^N
 \label{4.8}
 }
  where \alg{
   {\bf v}_k^N & = \bmat {cccc}
v_{kN+N-1}^T & v_{kN+N-2}^T & \ldots & v_{kN}^T \emat^T \nnl {\bf u}_k^N &=  \bmat {cccc}
u_{kN+N-1}^T & u_{kN+N-2}^T & \ldots & u_{kN}^T \emat^T \nnl
 {\bf
y}_k^N & =  \bmat {cccc} y_{kN+N-1}^T & y_{kN+N-2}^T & \ldots &
y_{kN}^T \emat^T .\nn} In model
(\ref{4.8}) we have \alg{ {\cal R}_N & =  \bmat {cccc} B & AB
& \ldots & A^{N-1}B \emat \nnl {\cal O}_N & = \bmat {cccc}
(CA^{N-1})^T & \ldots & (CA)^T & C^T \emat^T \nn \\
{\cal O}_N^R & = \bmat {cccc}
(A^{N-1})^T & \ldots & (A)^T & I
\emat^T  \nn\\
{\cal D}_N & = I_N \otimes D.
 \nn }
Note that, ${\cal R}_N$ and ${\cal O}_N$ denote,
respectively, the $N$-block reachability and observability matrices
of model  (\ref{state_space_model}), where the blocks forming
${\cal O}_N$ are written from bottom to top instead of the usual
top to bottom convention. In (\ref{4.8}), if
\alg{
H_k & =  \left \{ \begin{array} {cc}
CA^{k-1}B & k \geq 1 \\
0 & \mbox{otherwise }
\end{array} \right. \nn\\
  L_k & = \left \{ \begin{array} {cc}
A^{k-1}B & k \geq 1 \\
0 & \mbox{otherwise}
\end{array} \right.\nn
}
${\cal H}_N$ and ${\cal L}_N$ are block
{\em Toeplitz} matrices defined as follows
\alg{
{\cal H}_N & =  \bmat {cccccc}
0 & H_1 & H_2 & \cdots & H_{N-2} &H_{N-1}\\
0 & 0   & H_1 &  H_2 & \cdots & H_{N-2}\\
0 & 0   &  0 & H_1 &  \cdots & H_{N-3}\\
\vdots & \vdots &  \vdots &  & \vdots & \vdots\\
0 & 0   & \cdots & \cdots & 0 & H_1 \\
0 & 0   & \cdots & \cdots & \cdots & 0
\emat\nn\\
{\cal L}_N &= \bmat {cccccc}
0 & L_1 & L_2 & \cdots & L_{N-2} &L_{N-1}\\
0 & 0   & L_1 &  L_2 & \cdots & L_{N-2}\\
0 & 0   &  0 &  L_1 &  \cdots & L_{N-3}\\
\vdots & \vdots &  \vdots &  & \vdots & \vdots\\
0 & 0   & 0 & \cdots & \cdots & L_1 \\
0 & 0   & 0 & \cdots & \cdots & 0\\
\emat  \nn.}
We define \alg{{\cal J}_N &= {\cal O}_N^R  -{\cal L}_N{\cal H}_N^T [{\cal D}_N{\cal D}_N^T+{\cal H}_N{\cal
H}_N^T]^{-1}{\cal O}_N\nn\\
\Omega_N & = {\cal O}_N^T({\cal D}_N{\cal D}_N^T+{\cal H}_N{\cal H}_N^T)^{-1} {\cal O}_N\nn. }
Along similar lines used in \cite{ZORZI_CONTRACTION_CDC}, it is not difficult to see that the time-varying Riccati iteration associated to the
downsampled model (\ref{4.8}) takes the form
$P_{k+1}^d = r_{\tau,c,k}^d (P_k^d) $
where
\alg{ \label{downr}r_{\tau,c,k}^d (P_k^d) &:= \alpha_{N,k}
[ (P_k^d)^{-1} + \Omega_{\bar \Phi_{N,k}} ]^{-1}  \alpha_{N,k}^T
+ W_{\bar \Phi_{N,k}} \\
 \Omega_{\bar \Phi_{N,k}}
&= \Omega_N + {\cal J}_N^T  S_{\bar \Phi_{N,k}}^{-1} {\cal J}_N
\nn\\
 W_{\bar \Phi_{N,k}} &={\cal R}_N {\cal Q}_{\bar \Phi_{N,k}} {\cal R}_N^T\nn
}
with
\alg{ S_{\bar \Phi_{N,k}}& = - \bar \Phi_{N,k}^{-1}  + {\cal L}_N(I_{Nm}+{\cal H}_N^T ({\cal D}_N {\cal D}_N^T)^{-1}{\cal
H}_N)^{-1}{\cal L}_N^T \nn\\
{\cal Q}_{\bar \Phi_{N,k}} &= [ I_{Nm} + {\cal H}_N^T ({\cal D}_N {\cal D}_N^T)^{-1} {\cal H}_N -  {\cal L}_N^T\bar \Phi_{N,k}{\cal
L}_N ]^{-1}\nn\\
\bar  \Phi_{N,k} &  =  \mathrm{blockdiag}(\Phi_{kN+N-2} ,\Phi_{kN+N-3} ,\ldots , \Phi_{kN-1} ).\nn\\
\alpha_{N,k}&= A^N-\Rc_N(\Hc_N^T \Kc_{\bar  \Phi_{N,k}}^{-1} \Oc_N +\Lc_N^T  \Kc_{\bar  \Phi_{N,k}}^{-1} \Oc_N^R)\nn\\
 \Kc_{\bar  \Phi_{N,k}}&= \left[\begin{array}{cc} \Dc_N \Dc_N^T+\Hc_N \Hc_N^T & \Hc_N \Lc_N^T \\
 \Lc_N \Hc^T & -\bar\Phi_{N,k}^{-1}+\Lc_{N}\Lc_N^N \end{array}\right]
\nn} where we exploited the fact that $ \Dc_N   \Hc_N^T=0$ and $ \Dc_N   \Lc_N^T=0$ because $BD^T=0$.

\prop \label{prop_gramians} Let
\[ \tilde \phi_N=\frac{1}{\sigma_1({\cal L}_N(I_{Nm}+ {\cal H}_N^T ({\cal D}_N{\cal D}_N^T)^{-1}{\cal H}_N)^{-1}{\cal L}_N^T)}>0.\]
Assume that the pairs $(A,B)$ and $(A,C)$ are reachable and observable, respectively.
Then, there exits $\phi_N$, with $0 <  \phi_N < \tilde \phi_N$ and $N\geq n$, such that if $0\leq \bar \Phi \leq   \phi_N I_{nN}$
 then $\Omega_{\bar \Phi}$ and
$W_{\bar \Phi}$ are positive definite.
\eprop
\proof It is not difficult to see that $Q_{\bar \Phi}$ is positive definite and $S_{\bar \Phi}$ negative definite for $0\leq \bar \Phi < \tilde  \phi_N I_{nN}$. 
The mapping $\bar \Phi \mapsto W_{\bar \Phi}$ is nondecreasing with respect to the partial order of symmetric matrices over $0< \bar\Phi < \tilde \phi_n I_{nN}$
because its first variation along a direction $\delta\bar\Phi \in\bar\Qc_+^{Nn}$ is
\alg{\delta (W_{\bar \Phi};\delta\bar\Phi)= \Rc_N Q_{\bar \Phi} \Lc_N^T\delta \bar\Phi \Lc_NQ_{\bar \Phi}\Rc_N^T\geq 0.\nn}
 Note that,
$ W_{\bar \Phi=0} = {\cal R}_N ( I_{Nm}+{\cal H}_N^T ({\cal D}_N{\cal D}_N^T)^{-1}{\cal H}_N )^{-1}{\cal R}_N^T $
which is positive definite for $N\geq n$
because the pair $(A,B)$ is
reachable and thus $\Rc_N$ has full row rank. Accordingly, $W_{\bar \Phi}$ is positive definite for $0<\bar \Phi< \tilde \phi_N I_{nN}$. The mapping $\bar \Phi \mapsto \Omega_{\bar \Phi}$
is nonincreasing for $0<\bar \Phi<\tilde \phi_N I_{nN}$ because its first variation along $\delta\bar\Phi\in\bar\Qc_+^{Nn}$ is \alg{\delta(\Omega_{\bar\Phi};\delta\bar\Phi)=-\Jc_N^TS_{\bar \Phi}^{-1}  \bar\Phi^{-1} \delta \bar \Phi \bar\Phi^{-1} S_{\bar \Phi}^{-1}\Jc_N\leq 0.\nn} Moreover,
$  \Omega_{\bar \Phi=0}=\Omega_N$
which is positive definite for $N\geq n$
because the pair $(A,C)$ is observable. Accordingly, there exists a constant $\phi_N$ such that $0<\phi_N<\tilde \phi_N$
 and both $W_{\bar \Phi}$ and  $\Omega_{\bar \Phi}$ are positive definite for $0<\bar \Phi\leq \phi _N I_{Nn}$.\qed\\

\begin{remark} By the proof of Proposition \ref{prop_gramians}, one can see that $\phi_N$ can be computed as follows:
set $\phi_N=\tilde \phi_N$ and
check whether
$\Omega_{\phi_N I_{Nn}}$ is positive definite or not. If not, we decrease $\phi_N$ until $\Omega_{\phi_N I_{Nn}}$ becomes positive definite.\end{remark}

By Lemma \ref{lemma_lee_lim}, the mapping $r_{\tau,c,k}^d(\cdot)$ is strictly contractive provided that the matrices $ \Omega_{\bar \Phi_{N,k}} $
and $W_{\bar \Phi_{N,k}}$ are positive definite.  In view of Proposition \ref{prop_gramians}, if for some fixed $ \tilde q >0$ the following condition holds
\alg{\label{cond_bar_phi} \bar \Phi_{N,k}\leq\phi_N I_{nN}, \;\; k\geq \tilde q,} then the $N$-fold composition 
$r^d_{\tau,c,k}(\cdot)$ is strictly contractive for $k\geq \tilde q$ and thus $r_{\tau,c}(\cdot)$ is strictly contractive as well.

\section{Characterization of the range of the tolerance} \label{sec_toler}
 In this Section, we characterize a range of $c$
for which condition (\ref{cond_bar_phi}) holds.
The proofs of this Section only consider
the case $0<\tau<1$ because the results for the case $\tau=1$ can be proved along similar lines, and the case $\tau=0$ has been already proved in \cite{ZORZI_CONTRACTION_CDC}. 
Condition (\ref{cond_bar_phi}) is equivalent to the condition \alg{\label{cond_bar_phi2}\Phi_k\leq \phi_N I_{n},\;\; k\geq q+1}
for some $q>0$ fixed. Through  the next two Lemmas
we will  be able to derive a condition on $\theta_k$ which implies condition (\ref{cond_bar_phi2}).

\lem \label{prop_lower_bound} Let
$\bar P_{k+1}=r(\bar P_k)$, with $P_0=BB^T$,
be the convergent iteration generated by the usual Riccati mapping
\alg{r(P_k):=A (P_k^{-1}+ C^T(DD^T)^{-1}C)A^T+B B^T. \nn}
Consider the sequence generated by (\ref{iteration}).
Then, $P_k\geq \bar P_q$, with $k\geq q +1$
for any $q\geq 0$. \elem

\proof  It is well known that the sequence $\{\bar P_k\}$ is
 nondecreasing with respect to the partial order of the symmetric matrices. Accordingly, it is sufficient to prove that
$P_{k+1}\geq \bar P_k , \;\; k\geq 0$.
For this aim, we define the risk sensitive Riccati mapping, \cite{RISK_WHITTLE_1980},
\alg{r_{\Phi}^{RS}(P_k):=A (P_k^{-1}-\Phi +C ^T(DD^T)^{-1}C)A^T+B B^T \nn}
where $\Phi$ is a positive semidefinite matrix.
For $k=0$, we have $P_{1}=r_{\tau,c} (P_0)\geq B B^T =\bar P_0$.
Assume that $P_{k}\geq \bar P_{k-1}$, then
\alg{P_{k+1}&=r_{\tau,c}(P_k)=r^{RS}_{\Phi_{k-1}} (P_k)  \geq r(P_k)\geq r(\bar P_{k-1})=\bar P_k  \nn }
where we exploited the fact that $r^{RS}_{\Phi} (P)\geq r(P)$ for any $\Phi$ positive semidefinite and $P$ such that $0<P<\Phi^{-1}$, \cite{RISK_WHITTLE_1980},  and the fact that $r(\cdot)$ is a nondecreasing mapping with respect to the partial order of the symmetric matrices.  \qed\\

\lem\label{prop_upper_Phi_t}  Let $\bar d$ be such that
$P_{k+1}\geq \bar d I_n>0$, then
\alg{\Phi_k \leq  \left\{\begin{array}{cc}\frac{1-(1-\theta_k(1-\tau) \bar d)^{\frac{1}{1-\tau}}}{\bar d} I_n &  0\leq\tau<1 \\ \frac{1-\exp(-\theta_k\bar d)}{\bar d} I_n & \tau=1 . \\\end{array}\right. \nn }
\elem

\proof Consider the function \alg{f_\theta(\bar d)=\frac{1-(1-\theta(1-\tau) \bar d)^{\frac{1}{1-\tau}}}{\bar d}\nn }
defined over the set  $ \Sc=\{\bar d \hbox{ s.t }   0<\bar d< (\theta(1-\tau))^{-1}\}$ and $\theta>0$.
Then, \alg{\label{first_derivative_f}\frac{\mathrm{d}}{\mathrm{d}\bar d} f_\theta(\bar d) =\frac{-1+g_\theta(\bar d)}{\bar d^2}  }
where
\alg{g_\theta(\bar d)=(1-\theta(1-\tau)\bar d)^{\frac{1}{1-\tau}-1}(1+\theta \tau \bar d ). \nn}
It is not difficult to see that
\alg{\frac{\mathrm{d}}{\mathrm{d}\bar d} g_\theta(\bar d)=-\theta^2\bar d \tau (1-\theta(1-\tau) \bar d)^{\frac{1}{1-\tau}-2} \nn }
which  is nonpositive for $\bar d\in\Sc$. Accordingly, $g$ is a nonincreasing function over $\Sc$ and
\alg{ g_\theta(\bar d)\leq \underset{\bar d\rightarrow 0^+}{\lim} g_\theta(\bar d)= 1.\nn } Accordingly,  the first derivative of $f_\theta$ in (\ref{first_derivative_f})
is nonpositive over $\Sc$, i.e. $f_\theta$ is nonincreasing over $\Sc$.

 Let $L_{P_{k+1}}=\tilde U_{k+1} D_{k+1}^{\frac{1}{2}} U_{k+1}^T$ be the singular value decomposition of $L_{P_{k+1}}$, hence $\tilde U_{k+1}\tilde U_{k+1}^T=I_n$, $U_{k+1} U_{k+1}^T=I_n$
and $D_{k+1}^{\frac{1}{2}}= \mathrm{diag}(\ldots \,d_{i,k+1}^{\frac{1}{2}}\, \ldots)$ positive definite. Therefore, we have
\alg{ &V_{k+1}^{-1}= L_{P_{k+1}}^{-T} \left( I_n-\theta_{k} (1-\tau)U_{k+1} D_{k+1}U_{k+1}^T\right)^{\frac{1}{1-\tau}}L_{P_{k+1}}^{-1} \nn\\
&= L_{P_{k+1}}^{-T} \left( U_{k+1}U_{k+1}^T-\theta_{k} (1-\tau)U_{k+1} D_{k+1} U_{k+1}^T\right)^{\frac{1}{1-\tau}}L_{P_{k+1}}^{-1}\nn\\
&= L_{P_{k+1}}^{-T}  U_{k+1}\left( I_n-\theta_{k} (1-\tau) D_{k+1}\right)^{\frac{1}{1-\tau}} U_{k+1}^TL_{P_{k+1}}^{-1} \nn\\
&= \tilde U_{k+1} \mathrm{diag}\left(\ldots, \frac{(1-\theta_k(1-\tau) d_{i,k+1})^{\frac{1}{1-\tau}}}{d_{i,k+1}},\ldots\right) \tilde U_{k+1}^T.\nn }
Since the singular value decomposition of $P_{k+1}$  is  $P_{k+1}=\tilde U_{k+1}  \mathrm{diag} \left( \ldots, d_{i,k+1},\ldots\right)\tilde U_{k+1}^T$, we have \alg{\Phi_k &=P_{k+1}^{-1}-V_{k+1}^{-1}\nn\\
&=\tilde U_{k+1}  \mathrm{diag}\left(\ldots \, \frac{1- (1-\theta_k(1-\tau) d_{i,k+1})^{\frac{1}{1-\tau}}}{d_{i,k+1}} \,\ldots\right)\tilde U_{k+1}^T\nn\\
&= \tilde U_{k+1} \mathrm{diag}(\ldots,f_{\theta_k}(d_{i,k+1}),\ldots )  \tilde U_{k+1}^T. \nn }  By assumption,
$\bar d\leq d_{i,k+1}$, $i=1\ldots n$, therefore  we have
   $f_{\theta_k}(d_{i,k+1}) \leq f_{\theta_k} (\bar d)$, $i=1\ldots n$.
Accordingly, 
   $ \Phi_k \leq  \tilde U_{k+1}\mathrm{diag}(\ldots,f_{\theta_k}(\bar d),\ldots )  \tilde U_{k+1}^T=f_{\theta_k}(\bar d) I_n$
   which concludes the proof. \qed\\
   
Fixed $q>0$, by Lemma \ref{prop_lower_bound}, for the sequence generated by (\ref{iteration})
 we have $ P_{k} \geq \bar P_q\geq \sigma_n(\bar P_q) I_n$, $ \forall \, k\geq q+1$,  and by Lemma  \ref{prop_upper_Phi_t}
we have \alg{  \Phi_k  \leq \frac{1-(1-\theta_k(1-\tau)  \sigma_n(\bar P_q))^{\frac{1}{1-\tau}}}{\sigma_n(\bar P_q)} I_n,\;\;\; \forall \, k\geq q+1.\nn}
Therefore, the condition
\alg{ \frac{1-(1-\theta_k(1-\tau)  \sigma_n(\bar P_q))^{\frac{1}{1-\tau}}}{\sigma_n(\bar P_q)}  \leq \phi_N, \nn }
or equivalently
\alg{\label{relation_theta_t_phi_N}\theta_k \leq \frac{1-(1- \sigma_n(\bar P_q)\phi_N)^{1-\tau}}{(1-\tau) \sigma_n(\bar P_q)} , }
implies (\ref{cond_bar_phi2}). In particular, for $\tau=1$ we obtain
\alg{ \theta_k \leq \frac{-\log(1- \sigma_n(\bar P_q)\phi_N)}{ \sigma_n(\bar P_q)}.\nn}

The next Lemma is needed to derive a condition on $c$ which implies
condition (\ref{relation_theta_t_phi_N}), and thus also condition (\ref{cond_bar_phi2}).

\lem \label{prop_gamma}
Assuming that $0<\theta< ((1-\tau)\|P\|)^{-1}$, the following facts hold:\begin{enumerate}
                           \item $\gamma_\tau(P,\cdot)$ is monotone increasing over $\mathbb{R}_+$
                           \item If $ P\geq  Q$ then $\gamma_\tau(P,\theta)\geq \gamma_\tau(Q,\theta)$
                           \item $\gamma_\tau(P,\theta)>0$ for any $P\in\bar{\Qc}_n^+$ with $P\neq 0$.
                         \end{enumerate}
\elem

\proof  1) The statement has been proved in \cite{OPTIMALITY_ZORZI}.\\
2) First, note that
\alg{\gamma_\tau(P,\theta)=\tr \left(-\frac{1}{\tau(1-\tau)}(I_n-\theta(1-\tau) P)^{\frac{\tau}{\tau-1}} 
+\frac{1}{1-\tau} (I_n-\theta(1-\tau)P )^{\frac{1}{\tau-1}}+\frac{1}{\tau}I_n\right) . \nn }
To prove the statement, we show that the first variation
of $\gamma_\tau (P,\theta)$ with respect to $P$ in any direction $Q\in \bar \Qc_n^+$ is nonnegative:
\alg{ \delta \gamma_\tau  (P,\theta; Q)&=\frac{\theta}{1-\tau}\tr(-(I-\theta(1-\tau) P)^{\frac{1}{\tau-1}}Q
+(I-\theta(1-\tau) P)^{\frac{2-\tau}{\tau-1}} Q)\nn\\
& =\theta^2 \tr(P(I-\theta(1-\tau) P)^{\frac{2-\tau}{\tau-1}}Q )\nn\\
&=\theta^2 \tr(P^{\frac{1}{2}}(I-\theta(1-\tau) P)^{\frac{2-\tau}{2(\tau-1)}}Q (I-\theta(1-\tau) P)^{\frac{2-\tau}{2(\tau-1)}}P^{\frac{1}{2}})\geq 0 \nn }
where we exploited the fact that $(I-\theta(1-\tau) P)^{\frac{2-\tau}{\tau-1}}$ and $P$ commutes. \\
3)   $\gamma_\tau(P,\theta)$ is equal to the $\beta$-divergence between the covariance matrices $(I_n-\theta (1-\tau)P)^{\frac{1}{\tau-1}}$ and $I_n$, \cite{BETA}. Since
$(I_n-\theta (1-\tau) P)^{\frac{1}{\tau-1}}\neq I_n$, we get $\gamma_\tau(P,\theta)>0$.
\qed\\

\begin{figure}[h]
\begin{center}
\includegraphics[width=0.8\columnwidth]{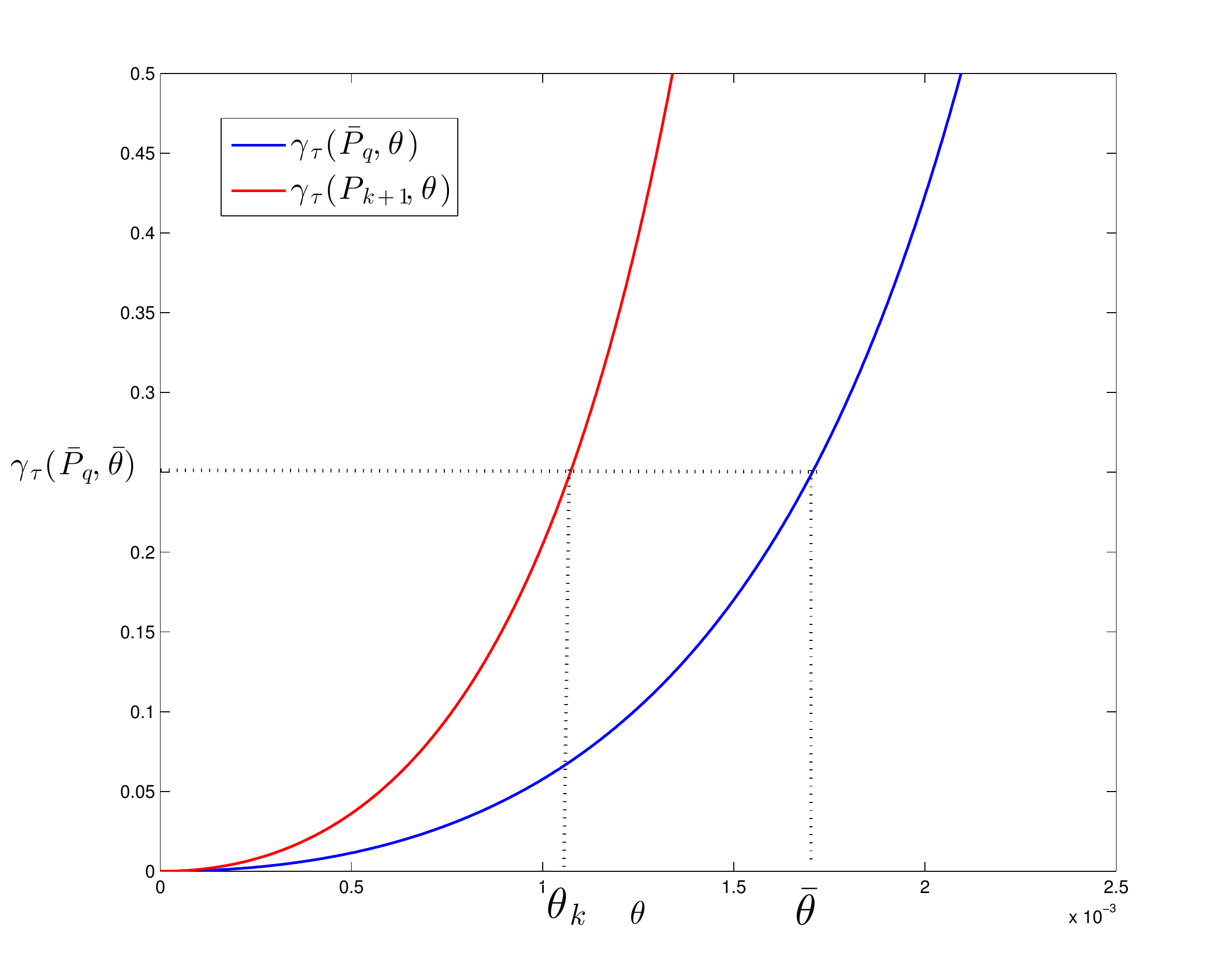}
\end{center}
 \caption{Function $\gamma_\tau$ with $P_{k+1}\geq \bar P_q$.}\label{gamma_proof}
\end{figure}

We know that $P_{k+1}\geq \bar P_q$ $\forall \, k\geq q$, which is equivalent to say $P_{k}\geq \bar P_q$ $\forall \, k\geq q+1$. Then, by Lemma \ref{prop_gamma}, condition $\gamma_\tau(P_{k+1},\theta_k)=\gamma_\tau( \bar P_q,\bar \theta)$ implies that
\alg{\theta_{k}\leq \bar \theta,\;\; \forall \, k\geq q. \nn }
Figure \ref{gamma_proof} shows this situation. Thus, (\ref{cond_bar_phi2}) holds if we choose $c$ in a such way that $\bar \theta\leq \frac{1-(1- \sigma_n(\bar P_q)\phi_N)^{1-\tau}}{(1-\tau) \sigma_n(\bar P_q)}$.

\teo \label{corollary} Let model (\ref{state_space_model}) be such that $(A,B)$ and $(A,C)$ are reachable and observable, respectively. Let $c$ be such that $0< c\leq c_{MAX}$ with \alg{c_{MAX} =
\left\{\begin{array}{cc} \gamma_\tau\left(\bar{P}_q, \frac{1-( 1-\sigma_n(\bar P_q)\phi_N)^{1-\tau}}{(1-\tau) \sigma_n(\bar P_q)}\right) & 0 \leq \tau<1 \\
\gamma_1\left(\bar{P}_q, \frac{-\log( 1-\sigma_n(\bar P_q)\phi_N)}{ \sigma_n(\bar P_q)}\right)  & \tau=1 \end{array}\right. \nn
} $N\geq n$ and $q>0$ are fixed. Then, for any $V_0\in\Qc_+^n$, the sequence $P_k$ generated by iteration (\ref{iteration}) converges to a unique solution $P$. Moreover, the limit $G$ of the filtering gain $G_k$ as $k\rightarrow \infty$ has the property that $A-GC$ is stable. 
\eteo

\proof Since \alg{c\leq \gamma_\tau\left(\bar{P}_q,\frac{1-(1-\sigma_n(\bar P_q)\phi_N)^{1-\tau}}{(1-\tau) \sigma_n(\bar P_q)}\right), \nn} by Lemma \ref{prop_gamma} we have that (\ref{relation_theta_t_phi_N}) holds  for $k\geq q$ and therefore
$\bar \Phi_{N,k}\leq \phi_N I_{nN}$  for $k\geq \tilde q=\lceil \frac{q}{N}\rceil$.
Accordingly, the mapping $r^d_{\tau,c,k}(\cdot)$ is strictly contractive for $k\geq \tilde q$.  Since $r^d_{\tau,c,k}(\cdot)$ is the $N$-fold composition of $r_{c,\tau}(\cdot)$, it follows that the sequence $P_k$ generated by (\ref{iteration}) converges. By (\ref{cond_on_theta_t}) the convergence of $P_k$ implies the convergence of $\theta_k$ to a unique value $\theta$. Thus, (\ref{def_V_t+1}) implies the convergence of  $V_k$ to a unique solution $V$. Finally, the stability of $A-GC$ can be proved by applying the Lyapunov stability theory to the algebraic Riccati-like equation 
\alg{P=(A-GC)V(A-GC)^T+BB^T+GG^T. \nn}
\qed\\

Finally, it is not difficult to show that the mapping \[ q\mapsto \gamma_\tau\left(\bar{P}_q,\frac{1-(1-\sigma_n(\bar P_q)\phi_N)^{1-\tau}}{(1-\tau)\sigma_n(\bar P_q)}\right) \]
is nondecreasing. Thus, we have to choose $q$ sufficiently large in order to find a bigger $c_{MAX}$.

\section{Example} \label{sec_example}
\begin{figure}[h]
\begin{center}
\includegraphics[width=0.8\columnwidth]{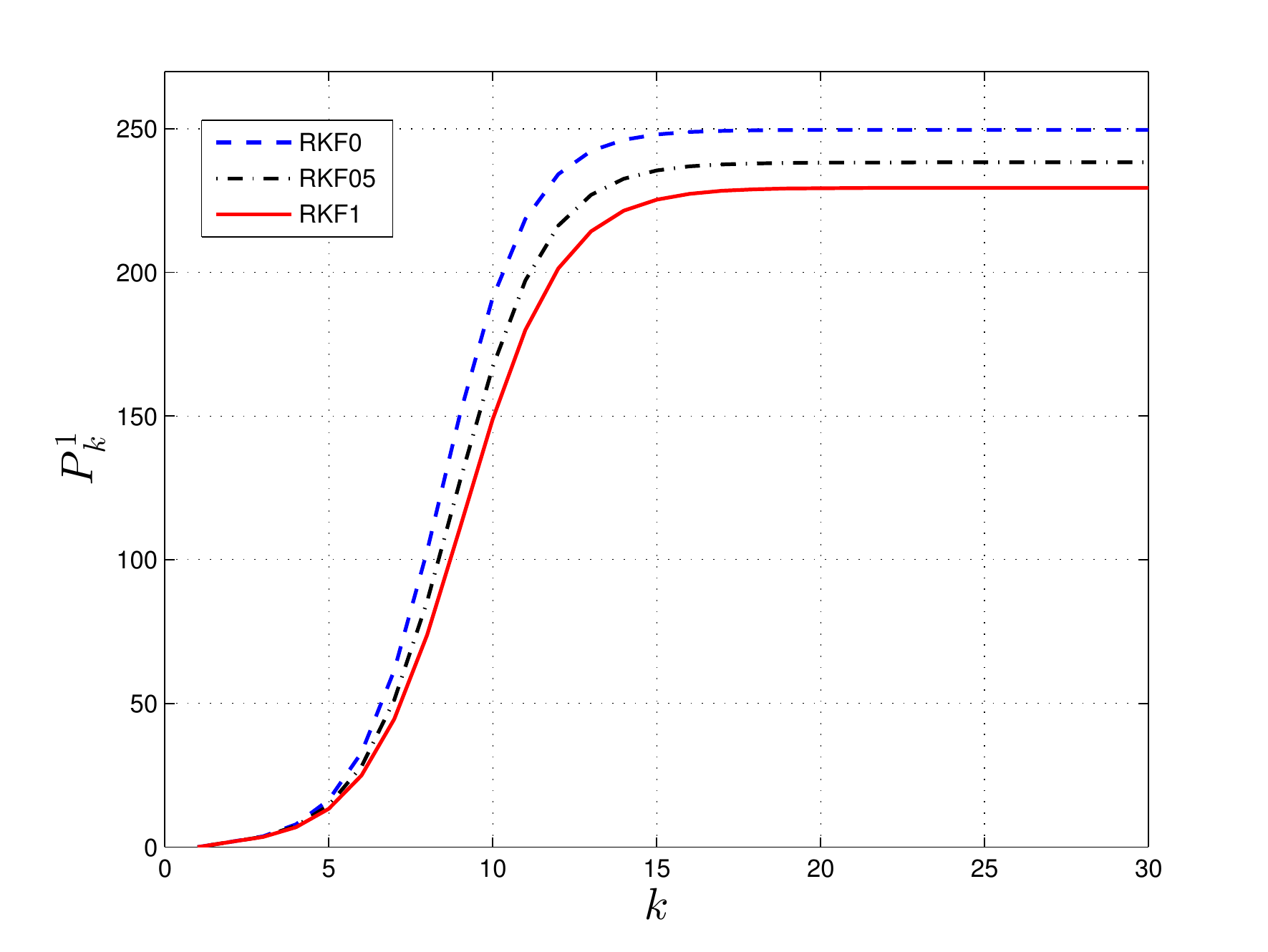}
\end{center}
 \caption{Pseudo-nominal variance of the state estimation error  of the first component of $x_k$.}\label{nomiP1}
\end{figure}

\begin{figure}[h]
\begin{center}
\includegraphics[width=0.8\columnwidth]{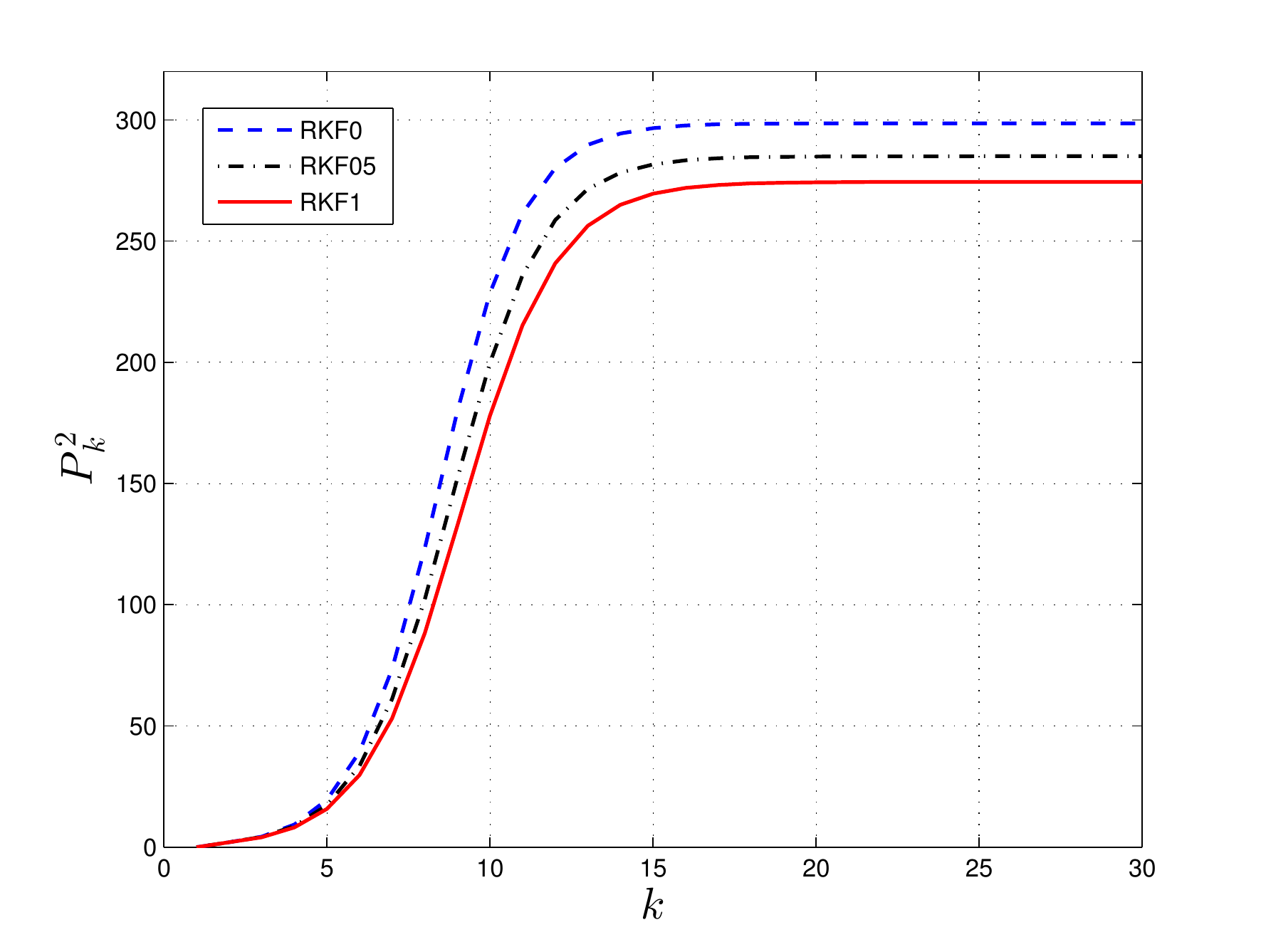}
\end{center}
 \caption{Pseudo-nominal variance of the state estimation error of the second component of $x_k$.}\label{nomiP2}
\end{figure}

\begin{figure}[h]
\begin{center}
\includegraphics[width=0.8\columnwidth]{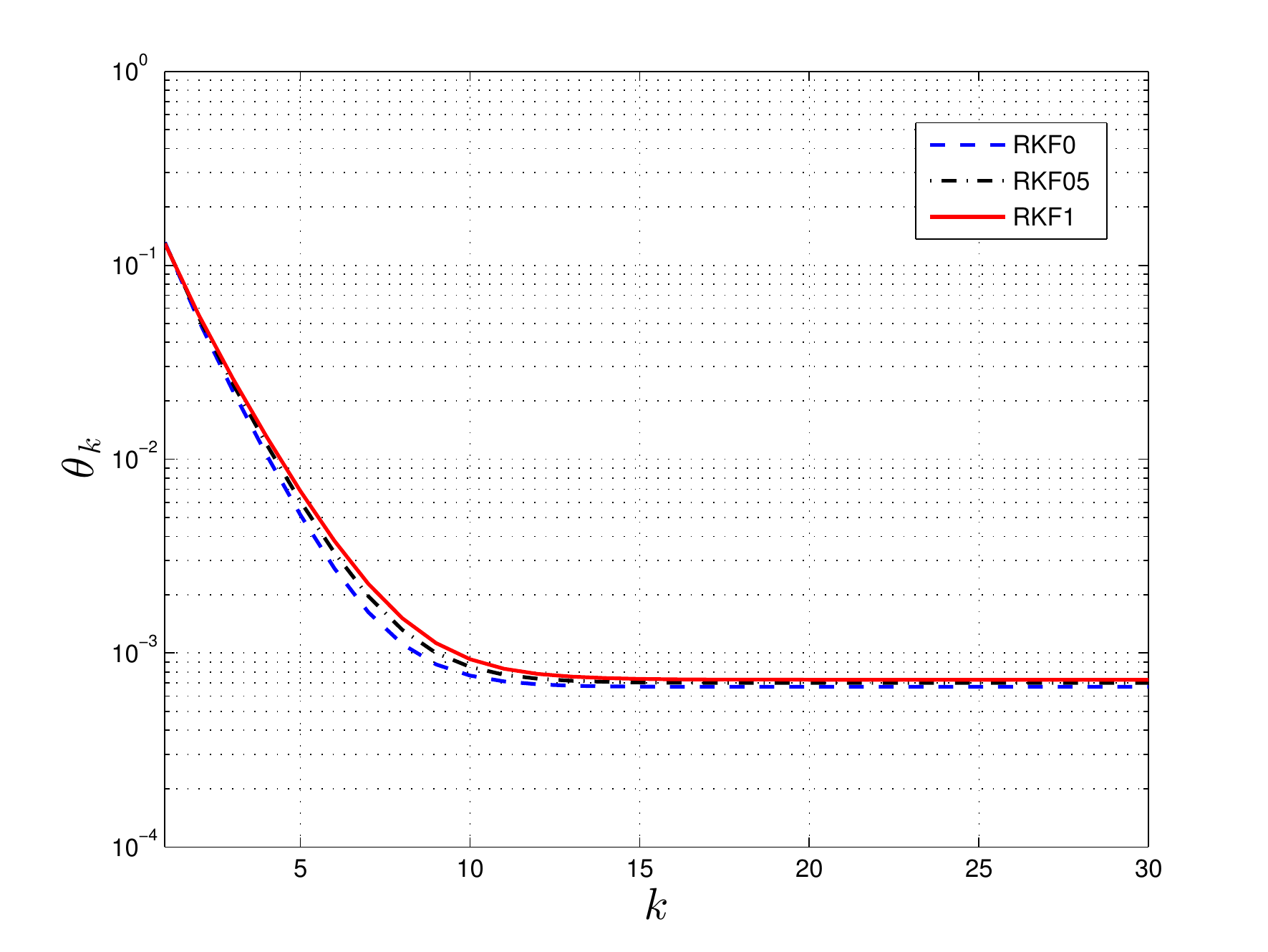}
\end{center}
 \caption{Time-varying risk-sensitivity parameter $\theta_k$ in logarithmic scale.}\label{theta}
\end{figure}

\begin{figure}[h]
\begin{center}
\includegraphics[width=0.8\columnwidth]{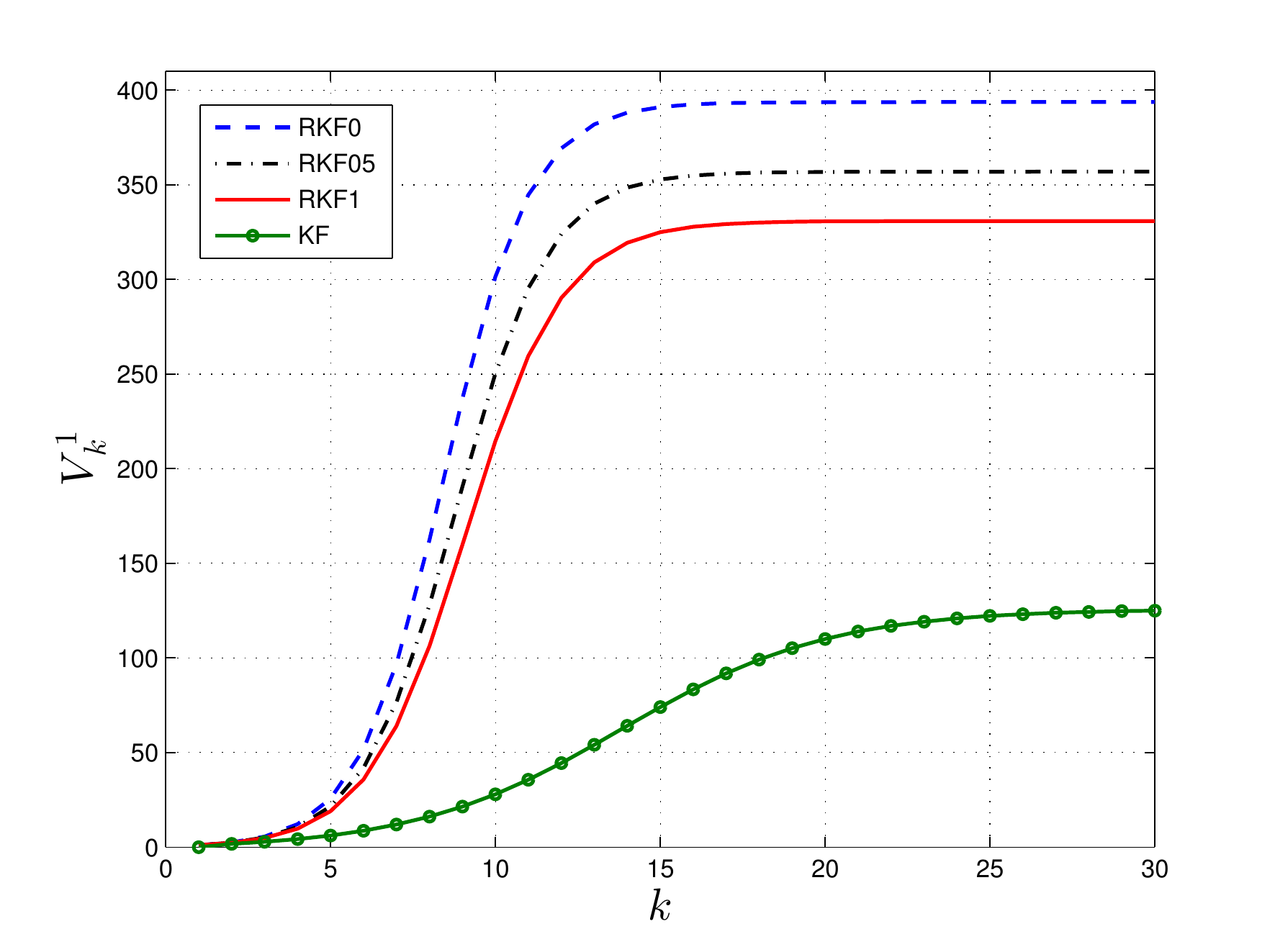}
\end{center}
 \caption{Least-favorable variance of the state estimation error  of the first component of $x_k$.}\label{nomiV1}
\end{figure}

\begin{figure}[h]
\begin{center}
\includegraphics[width=0.8\columnwidth]{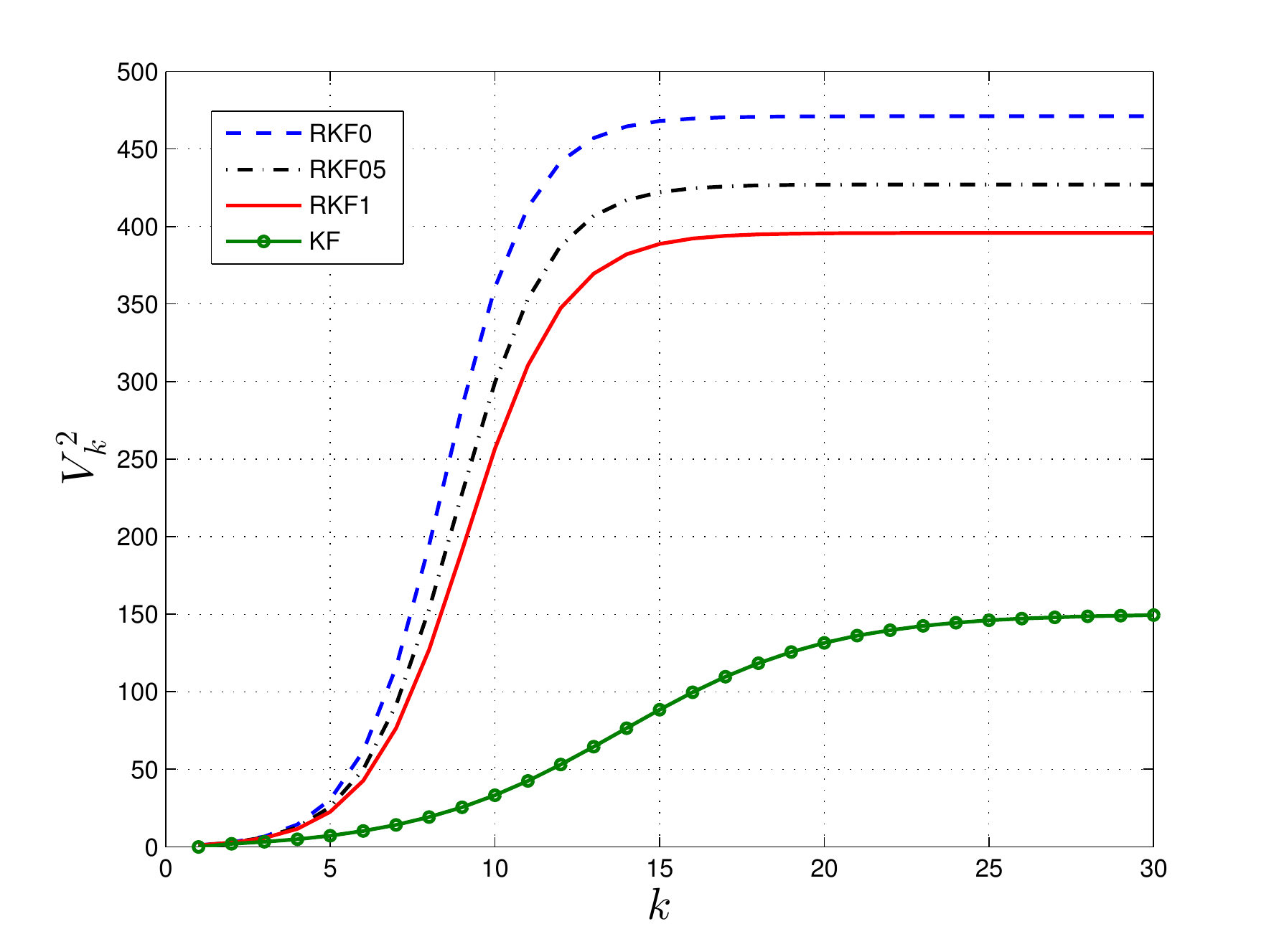}
\end{center}
 \caption{Least-favorable variance of the state estimation error of the second component of $x_k$.}\label{nomiV2}
\end{figure}

We consider the constant state space model (\ref{state_space_model}) used in \cite{LEVY_ZORZI_RISK_CONTRACTION},
\alg{ A&=\left[\begin{array}{cc}0.1 & 1  \\0 & 1.2 \\\end{array}\right], \;\; B=\left[\begin{array}{ccc} 1 & 0 & 0  \\0 & 1 &0 \\\end{array}\right]\nn\\
C&=\left[\begin{array}{cc} 1 & -1  \\ \end{array}\right], \;\; D=\left[\begin{array}{ccc} 0 & 0  & 1  \\ \end{array}\right]. \nn}
The error covariance matrix at time $k=0$ is chosen as $V_0=I_2$. 
We study the convergence of filter (\ref{rob_est_formula1})-(\ref{def_gamma}) with three different values for $\tau$: $\tau=0$, $\tau=0.5$ and $\tau=1$. Fixing $q=40$, $N=50$
we found that
\alg{
\bar P_q =10^2\cdot \left[\begin{array}{cc} 1.2568 &  1.3641\\
   1.3641  & 1.5025\\\end{array}\right],\;\; \tilde \phi_N=
  1.3335\cdot 10^{-3},\;\;
\phi_N =   1.3328 \cdot10^{-3}. \nn} Moreover,
the robust Kalman filter (\ref{rob_est_formula1})-(\ref{def_gamma}) converges with tolerance in the range $[0,c_{MAX}]$ where
\alg{c_{MAX}&= 1.22\cdot 10^{-1}   \hbox{  for } \tau=0\nn\\
c_{MAX}&= 1.01\cdot 10^{-1}   \hbox{  for } \tau=0.5\nn\\
  c_{MAX} &= 8.62\cdot 10^{-2}  \hbox{  for } \tau=1. \nn}
Now, we compare the performances of the following three filters:
\begin{itemize}
\item KF: the standard Kalman filter
\item RKF0: the robust Kalman filter with  $\tau=0$ and $c=1.22\cdot 10^{-1}$
\item RKF05: the robust Kalman filter  with  $\tau=0.5$ and $c=1.01\cdot 10^{-1} $
\item RKF1: the robust Kalman filter with  $\tau=1$ and $c= 8.62\cdot 10^{-2}$
\end{itemize}
that is we consider the robust Kalman filter with $\tau=0$, $\tau=0.5$, $\tau=1$ with the corresponding maximum tolerance for which we know that it converges. In Figure \ref{nomiP1} we show the pseudo-nominal variance  
of the state estimation error of the first component of the state, that is the entry of $P_k$ in position (1,1). 
In Figure \ref{nomiP2} we show the pseudo-nominal variance of the state estimation error of the second component of $x_k$, that is the entry in position (2,2) of $P_k$.
 Roughly speaking these quantities represent the error variance computed using the nominal density $f_k$ but propagating the previous least favorable density $\tilde f_{k-1}$. 
The previous figures show that the Riccati-like iteration converges after 20 steps for $\tau=0$, $\tau=0.5$ and $\tau=1$. 
In Figure \ref{theta}, we show the time-varying risk-sensitivity parameter $\theta_k$ which after 20 steps is already constant.   
In Figure \ref{nomiV1}
and Figure \ref{nomiV2} we consider the corresponding least-favorable error variance, i.e. the error variance is computed by using the least-favorable density $\tilde f_k$ and propagating the previous least favorable 
density $\tilde f_{k-1}$.  It is clear that RKF0, RKF05 and RKF1 are very conservative with respect to the KF, i.e. their error variances are larger than the ones given by KF. This means that, although the upper bound $c_{MAX} $ we found is not tight, the range $[0,c_{MAX}]$ contains a sufficiently large class of robust estimators. In other words,  with c close to zero we have robust Kalman filters with performance similar to KF, while with $c$ close to $c_{MAX}$ we have robust Kalman filters very different than KF.

\section{Convergence analysis of the $\tau$-risk sensitive filters} \label{sec_RS}
Consider the state-space model (\ref{state_space_model}) and the corresponding nominal joint Gaussian probability density ${f}_k(x_{k+1},y_k|Y_{k-1})$. The family of risk sensitive filters \cite{STATETAU_2017} parametrized by $\tau\in[0,1]$ is given by
 \alg{\label{RS_est_minimax}\hat x_{k+1}=\underset{g_k\in\Gc_k}{\mathrm{argmin}}  \underset{\tilde f_k\in\Bc_{k,\tau}}{\max}\Es_{\tilde f_k}[ & \|x_{k+1}  -g_k(y_k)\|^2 \,|\,Y_{k-1}]-\theta^{-1}\Dc_{\tau}(\tilde f_k\| f_k)}
where  $\tilde f_k$ is Gaussian, $\Bc_{k,\tau}=\{\tilde f_k \hbox{ s.t. } \Ds_{\tau}(\tilde f_k \| f_k)<\infty  \}$ and  $\Gc_k$ is the set of estimators for which the objective function in (\ref{RS_est_minimax}) is finite. $\theta>0$ is the risk sensitivity parameter.  The second term in the objective function in (\ref{RS_est_minimax}) is always nonpositive because $\Ds_{\tau}(\tilde f_k\| f_k)\geq 0$.
Therefore, for large values of $\theta$ the maximizer has the possibility to take a probability density far from the nominal one. The $\tau$-risk sensitive filter (\ref{RS_est_minimax}) thus represents a relaxed version of the robust Kalman filter (\ref{rob_est_formula1})-(\ref{def_gamma}) where $\theta$ now is constant and fixed by the user. For the case $\tau=0$ we obtain the usual risk sensitive filter \cite{boel2002robustness}. The resulting estimator obeys the recursion (\ref{rob_est_formula1})-(\ref{def_P_t+1}) with 
 \alg{ & V_{k+1}=\left\{
      \begin{array}{ll}
        L_{P_{k+1}}\left(I_n-\theta(1-\tau) L_{P_{k+1}}^TL_{P_{k+1}}\right)^{\frac{1}{\tau-1}} L_{P_{k+1}}^T , & \hspace{-0.1cm}  0<\tau<1\\
L_{P_{k+1}}\exp\left(\theta L_{P_{k+1}}^T L_{P_{k+1}}\right) L_{P_{k+1}}^T, & \tau=1.
      \end{array}
    \right.\nn} 
The study of the asymptotic behavior of the $\tau$-risk sensitive filter requires to consider two different cases:  the case $0<\tau<1$  and the case $\tau=1$.    
    
 In the former case, the Riccati-like iteration has the same form of (\ref{iteration}) but the image of $\Qc_+^n$ under this mapping is not entirely contained in $\Qc_+^n$. The reason is that condition $V_{k}>0$ holds only if $P_k$ is such that $0<P_{k}<(\theta(1-\tau))^{-1}I_n$ and this condition could be not satisfied. Following similar arguments used in \cite{LEVY_ZORZI_RISK_CONTRACTION} for the case $\tau=0$, it is possible to find conditions on $V_0$ and $\theta$ for which the trajectory of iteration (\ref{iteration}) satisfies $V_k>0$ for any $k> 0$. However, these conditions on $V_0$ and $\theta$ are rather intricate and require to design a gain matrix and a scaling factor $\rho^2$.

For the case $\tau=1$, $V_k$ is positive definite, and thus well defined, whenever $P_{k}$ is positive definite. Accordingly, the image of $\Qc_+^n$ under the corresponding mapping, denoted by $r_{\tau=1,\theta}(\cdot)$, is $\Qc_+^n$.
Thus, the convergence of the iteration is guaranteed by only imposing conditions on the risk sensitivity  parameter $\theta$. 
\teo  Let model (\ref{state_space_model}) be such that $(A,B)$ and $(A,C)$ are reachable and observable, respectively. Let $\theta$ be such that 
\alg{\label{condthrs}\theta\leq \frac{-\log(1- \sigma_n(\bar P_q)\phi_N)}{ \sigma_n(\bar P_q)},  }
 $N\geq n$ and $q>0$ are fixed. Then, for any $V_0\in\Qc_+^n$, the sequence $P_k$ generated by 
the risk-sensitive filter with $\tau=1$ converges to a unique solution $P$. Moreover, the limit $G$ of the filtering gain $G_k$ as $k\rightarrow \infty$ has the property that $A-GC$ is stable. 
\eteo
\proof We consider the downsampled process $x_k^d$ with $x_k^d=x_{kN}$ and the corresponding time-varying Riccati-iteration is $P_{k+1}^d=r^d_{\tau,\theta,k}(P_k^d)$ where $r^d_{\tau,\theta,k}(\cdot)$ has the same structure of (\ref{downr}). Let $\Phi_k=P_{k+1}^{-1}-V_{k+1}^{-1}$. Proposition \ref{prop_gramians} still holds. In particular, there exists $\phi_N$ such that if  (\ref{cond_bar_phi2}) holds then the matrices $\Omega_{\bar \Phi_{N,k}}$ and $W_{\bar \Phi_{N,k}}$
are positive definite. Accordingly, by Lemma \ref{lemma_lee_lim} the $N$-fold mapping  $r^d_{\tau,\theta,k}(\cdot)$ is strictly contractive and thus also $r_{\tau=1,\theta}(\cdot)$ is strictly contractive. Lemma \ref{prop_lower_bound} and Lemma \ref{prop_upper_Phi_t} still hold, in particular \alg{\Phi_k\leq \frac{1-e^{-\sigma_n(\bar P_q)\theta}}{\sigma_n(\bar P_q)} I_n, \;\; \forall \; k\geq q+1.\nn} Finally, by imposing \alg{\frac{1-e^{-\sigma_n(\bar P_q)\theta}}{\sigma_n(\bar P_q)} \leq \phi_N,\nn} which coincides with (\ref{condthrs}), then condition (\ref{cond_bar_phi2}) holds. Thus, the sequence $P_k$ converges to a unique $P$ as $k\rightarrow \infty$. The stability of $A-GC$ 
follows as before.\qed\\

It is clear that condition (\ref{condthrs}) on the risk sensitivity parameter is easy to check. Accordingly, this filter is preferable than the risk-sensitive filter with $0\leq \tau<1$.

\section{Conclusions} \label{sec_conclusions}
A convergence analysis of a family of robust Kalman filters has been presented. This analysis exploited the fact that the $N$-fold Riccati mapping, which is given by downsampling these filters,  is strictly contractive provided that the time-varying risk-sensitive parameter is sufficiently small. This condition is then guaranteed by placing an upper bound on the tolerance parameter of the robust filters. Finally, we have studied the convergence property of a family of risk-sensitive filters which can be understood as a relaxed version of the previous robust Kalman filters.

\end{document}